\documentclass[12pt]{article}
\usepackage[all]{xy}
\usepackage{amssymb,amsmath,amsthm,amscd,amsfonts}
\usepackage{hyperref}
\newtheorem{df}{Definition}[section]
\newtheorem{thm}[df]{Theorem}
\newtheorem{pro}[df]{Proposition}

\newtheorem{lem}[df]{Lemma}
\newtheorem{cor}[df]{Corollary}
\newtheorem{rk}[df]{Remark}
\newtheorem{ex}[df]{Example}

\newcommand{\E}{{\cal E}}

\def\id{{\rm id}}

\date{}
\title{\bf  Lawvere-Tierney sheaves, factorization systems, sections and $j$-essential monomorphisms in a topos}
\author{{\bf Zeinab Khanjanzadeh}\and {\bf Ali Madanshekaf } \\
Department of Mathematics\\Faculty of Mathematics, Statistics and Computer Science\\
Semnan University\\
Semnan\\
Iran\\ emails:
z.khanjanzadeh@gmail.com\\ \qquad\qquad amadanshekaf@semnan.ac.ir}
\date{}
\begin{document}
\maketitle
\begin{abstract}
 Let  $j$ be  a Lawvere-Tierney topology  (a topology, for short)  on an arbitrary topos $\mathcal{E}$, $B$
 an object  of $\mathcal{E}$, and  $j_B =  j\times 1_B$  the induced topology  on
 the slice topos $\mathcal{E}/B$. In this manuscript, we analyze some properties of the  pullback functor
 $\Pi_B:\mathcal{E}\rightarrow \mathcal{E}/B$ which have deal with
 topology.
Then
 for a left cancelable class $\mathcal{M}$ of all
$j$-dense  monomorphisms in a
  topos $\mathcal{E}$, we  achieve some
   necessary and sufficient conditions for that $(\mathcal{M} , \mathcal{M}^{\perp})$ is a factorization system in  $\mathcal{E}$, which is related to  the
    factorization systems in  slice topoi $\mathcal{E}/B,$ where $B$ ranges over the class of objects of $\mathcal{E}$.
   Among other things, we  prove  that an arrow $f : X\rightarrow B$ in $\mathcal{E}$ is a
    $j_B$-sheaf whenever the graph of $f$, is a section in $\mathcal{E}/B$ as well as the object of sections $S(f)$ of $f$,  is a $j$-sheaf in $\mathcal{E}$.
    Furthermore,   we introduce  a class of monomorphisms in $\mathcal{E}$, which we call them $j$-essential.
 Some equivalent forms of those and some of their properties are  presented. Also, we prove that any presheaf in a presheaf topos
 has a maximal essential extension. Finally, some similarities and
 differences of the obtained result are discussed if we  put a  (productive)
 weak topology $j$, studied by some authors,  instead of a  topology.
\end{abstract}
AMS {\it subject classification}: 18B25; 18A25; 18A32; 18F20; 18A20. \\
{\it key words}:    (Weak) Lawvere-Tierney topology;  Sheaf;
Factorization system; Slice topos; Essential monomorphism.
\section{Introduction and background}
 A Lawvere-Tierney topology is a  logical connective for modal
logic. Recently, applications of Lawvere-Tierney topologies in
broad topics such as measure theory~\cite{measure} and quantum
Physics~\cite{quantum1,quantum2} are observed.  In the spacial
case, considerable work has
 been presented that is dedicated to the study  of (weak)  Lawvere-Tierney topology on a presheaf topos on a small category and especially on
 a monoid, see~\cite{hosseini,  Luis Espa}.
It is clear that Lawvere-Tierney sheaves in a topos  are exactly
injective objects  (of course, with respect to dense
monomorphisms, not to merely monomorphisms) which are separated
too. Injectivity with respect to a class $\mathcal{M}$ of
morphisms in a slice category $\mathcal{C}/B$ (which its objects
are $\mathcal{C}$-arrows with codomain $B$) has been
 studied in extensive form, for example we refer the reader to~ \cite{Adamek, Cagliari}. From this perspective, in this paper we will establish
 some  categorical characterizations of injectives in slice topoi to  sheaves.  The object of sections $S(f)$ of  $f$  is a notion which
 in~\cite{Cagliari} it is related to injective objects in a slice category. This object  is very useful in synthetic differential geometry (or SDG, for short) (for
details, see \cite{Lavendhome}). For example, considering $D$ as
infinitesimals,  for any micro-linear object  $M$ we have:
\begin{itemize}
\item Let $\tau$ be the tangent bundle on $M$, i.e., $\tau : M^D\rightarrow M$, which is defined by $\tau (t) = t(0).$
Then $S(\tau)$ is all vector fields on $M$.
\item Consider $\eta:M^{D\times D}\rightarrow M$ which assigns  to any micro-square $Q$ of $M^{D\times D},$
the element  $Q(0,0)$. Then, $S(\eta)$ is all distributions of dimension 2 on $M$.
\end{itemize}

 Throughout this paper, $\mathcal{E}$ is a (elementary) topos, two objects
$0, 1$ are the initial and terminal objects and  the object $\Omega$ together with the arrow $1\stackrel{{\rm true}}{\rightarrowtail} \Omega$ is the
subobject classifier  of $\mathcal{E}$. Also, the arrow $\wedge:\Omega\times\Omega\rightarrow \Omega$ is the meet operation on $\Omega$.  Now,
 we  express  some  basic concepts from~\cite{maclane} which will be needed in sequel.
\begin{df}\label{loc. ope.}
{\rm A {\it Lawvere-Tierney topology }  on $\mathcal{E}$ is a map  $j:\Omega\rightarrow \Omega$ in   $\mathcal{E}$ satisfies the following properties
\[\begin{array}{ccc}
  ({\rm a})~ j\circ {\rm true} = {\rm true}; &({\rm b})~ j\circ j= j; &({\rm c})~ j\circ \wedge =  \wedge\circ (j\times j); \\
\xymatrix{
1 \ar[dr]_{{\rm true}} \ar[r]^{\rm true} & \Omega \ar[d]^j \\
&\Omega } & \xymatrix{
\Omega \ar[dr]_j \ar[r]^j & \Omega \ar[d]^j \\
&\Omega } &\xymatrix{
\Omega\times\Omega\ar[r]^{~\wedge}  \ar[d]_{j\times j}& \Omega \ar[d]^{j}
 \\
\Omega\times\Omega \ar[r]^{~\wedge} & \Omega }
 \end{array}
\]
Form now on,  we say briefly to a Lawvere-Tierney topology on
$\mathcal{E}$, a {\it topology} on $\mathcal{E}$.}
\end{df}
Recall \cite{maclane}   that  topologies on $\mathcal{E}$  are in one to one correspondence with universal closure operators. For a topology
$j$ on $\mathcal{E}$, considering $\overline{(\,\cdot\,)}$ as the universal closure operator corresponding to $j$, a monomorphism $k:A\rightarrowtail C$
 in $\mathcal{E}$ is called {\it $j$-dense} whenever $\overline{A} = C$, as two subobjects of $C$. Also, we say that $k$ is {\it $j$-closed}
 if we have $\overline{A} = A$, again as subobjects of $C$.
\begin{df}\label{she.}
{\rm For a topology $j$ on $\mathcal{E}$, an object $F$ of $\mathcal{E}$ is called {\it a $j$-sheaf} whenever for any $j$-dense monomorphism
 $m:A\rightarrowtail E$, one can uniquely  extend any arrow $h:A\rightarrow F$  to a map $g$ on all of $E$,
\begin{equation}\label{ne.}
\SelectTips{cm}{}\xymatrix{ A \ar[r]^h\ar@{>->}[d]_{m} & F \\  E \ar@{-->}[ur]_{g}}
\end{equation}
We say that $F$ is {\it  $j$-separated} if the arrow $g$ exists in (\ref{ne.}), it is unique.}
\end{df}
We will denote the full subcategories  of $\mathcal{E}$ consisting of $j$-sheaves and  $j$-separated objects as  ${\bf Sh}_j(\mathcal{E})$
and  ${\bf Sep}_j(\mathcal{E})$, respectively.

We now briefly describe the contents of  other sections. We start
in Section 2, to study  basic properties of  the  pullback functor
 $\Pi_B:\mathcal{E}\rightarrow \mathcal{E}/B$, for any object  $B$  of $\mathcal{E}$, along with the unique map $!_B: B\rightarrow 1$.
 Afterwards, we  would like to achieve, for a left cancelable class $\mathcal{M}$ of all $j$-dense  monomorphisms in a topos $\mathcal{E}$,
 some necessary and sufficient conditions for that $(\mathcal{M} , \mathcal{M}^{\perp})$ to be a factorization system in  $\mathcal{E}$, which
  is related to  the  factorization systems in  slice topoi $\mathcal{E}/B$.
 In section 3, among other things, we  prove  that an arrow $f : X\rightarrow B$ in $\mathcal{E}$ is a $j_B$-sheaf
  whenever the graph of $f$, is a section in $\mathcal{E}/B$ as well as the object of sections $S(f)$ of $f$,  is a $j$-sheaf in $\mathcal{E}$.
   In section 4,  we introduce  a class of monomorphisms in an elementary topos $\mathcal{E}$, which we call them `$j$-essential monomorphisms'. We present
    some equivalent forms of these and some of their properties. Meanwhile, we prove that any presheaf in a presheaf topos has a maximal essential extension.
    It is shown that the functor $\Pi_B$ reflects $j$-essential extensions.   It is seen  that some of these results hold for a
     {\it (productive) weak topology} $j$, studied in~\cite{Madanshekaf}, instead of a  topology as well.
\section{Pullback functors, left cancelable dense monomorphisms and factorization systems}
The purpose of this section is to present some basic properties of
the  pullback functor  $\Pi_B:\mathcal{E}\rightarrow
\mathcal{E}/B$, for any object  $B$  of $\mathcal{E}$,  along with
the unique map $!_B: B\rightarrow 1$. Afterwards, for a left
cancelable class $\mathcal{M}$ of all $j$-dense  monomorphisms in
a topos $\mathcal{E}$ we  achieve some necessary and sufficient
conditions for that $(\mathcal{M} , \mathcal{M}^{\perp})$ to be a
factorization system in  $\mathcal{E}$, which is related to  the
factorization systems in  slice topoi $\mathcal{E}/B$.

To begin with, the following lemma  characterizes   sheaves   in
a topos $\mathcal{E}$.
\begin{lem}\label{unique retract}
Let  $j$ be a topology on $\mathcal{E}$. Then an object $E$ of
$\mathcal{E}$ is  $j$-sheaf iff  $E$ is {\it $j$-unique absolute
retract}; that is,  any $j$-dense monomorphism $u:E\rightarrowtail
F,$ has a unique retraction $v : F\rightarrow E$.
\end{lem}
{\bf  Proof.} {\it Necessity.} Since $E$ is a $j$-sheaf,   for any
$j$-dense monomorphism $u : E\rightarrowtail F$, corresponding to
the identity map $\id_E:E\rightarrow E$ there exists a unique map
$v : F\rightarrowtail E$ such that the following diagram commutes.
$$\SelectTips{cm}{}\xymatrix{E \ar@{>->}[d]_{u}\ar[r]^{\id_E} & E  \\ \ar@{-->}[ur]_{v}F}$$

{\it Sufficiency.} For each  $j$-dense monomorphism $m:U\rightarrowtail V$ and any map $f:U\rightarrow E$, we construct the following pushout
 diagram in $\mathcal{E}$.
\begin{equation}\label{12}
\SelectTips{cm}{}\xymatrix{   U  \ar@{>->}[d] _{m} \ar[r]^{f}& E\ar@{>->}[d]^{n} \\ V \ar[r]_{g~~}^{{~~~~{\rm p.o.}}} & F  }
\end{equation}
Since in any topos pushouts transfer $j$-dense monomorphisms (see
\cite{jhonstone}), so, in (\ref{12}), $n$ is $j$-dense and hence by
 assumption, there exists a unique retraction $p : F\rightarrow
E$ such that $pn = \id_E$. Now, for the the arrow $pg :
V\rightarrow E$ we have $pg m= pnf= \id_Ef=f.$ To prove that $pg :
V\rightarrow E$ with this  property  is unique, let $h :
V\rightarrow E$ be an arrow in $\mathcal{E}$ in such a way that
$hm = f$. Then, in the pushout diagram (\ref{12}), according to
the maps $h:V\rightarrow E$ and $\id_E:E\rightarrow E$,  there
exists a unique map $k : F\rightarrow E $ such that $kn=\id_E$ and
$kg = h$.
$$\SelectTips{cm}{}\xymatrix{ U \ar[r]^{f} \ar@{>->}[d]_{m} & E \ar@{>->}[d]^{n}\ar@/^1pc/[ddr]^{\id_E} \\
V \ar[r]_{g}\ar@/_1pc/[drr]_{h} & F\ar@{-->}[dr]^{k}\\ && E}$$
Now, $k$ is a retraction of
$j$-dense monomorphism $n,$ so by hypothesis  we get $p=k$.
Consequently, $pg = kg =h $. $\qquad\square$

 For an  object $B$ of $\mathcal{E}$, we consider the pullback functor  $\Pi_B:\mathcal{E}\rightarrow \mathcal{E}/B$ along with the
  unique map $!_B: B\rightarrow 1$, which assigns  to any $A$ of $\mathcal{E}$, the second projection $\Pi_B(A)=\pi_B^A : A\times B\rightarrow B$
  and to any $f : A\rightarrow C$, the arrow $f\times \id_B:A\times B\rightarrow C\times B$ in $\mathcal{E}$ such that $\pi_B^C (f\times \id_B) = \pi_B^A$.
   Recall \cite{maclane} that   the object $\pi_B^{\Omega}$ together with the arrow
 $${\rm true} \times \id_B:\id_B\longrightarrow \pi_B^{\Omega}$$ is the subobject classifier of the slice topos $\mathcal{E}/B$. Also, in a similar vein,
 we can observe that the meet operation $\wedge_B$ on $\pi_B^{\Omega}$ is the arrow
 $\wedge\times 1_B$ in $\mathcal{E}$ such that $\pi_B^{\Omega} (\wedge\times 1_B) = \pi_B^{\Omega\times\Omega},$
$$\xymatrix{\Omega\times\Omega\times B\ar[r]^{~~~\wedge\times 1_B}\ar[dr]_{\pi_B^{\Omega\times\Omega}} &\Omega\times B
\ar[d]^{\pi_B^{\Omega}}\\ &B.}$$
 Now, by Definition \ref{loc. ope.}, we easily get the following lemma.
\begin{lem}\label{lop(e/b)}
 Let $B$ be any  object in a topos $\mathcal{E}$. Then any topology $k:\pi_B^{\Omega}\rightarrow \pi_B^{\Omega}$ on
 $\mathcal{E}/B$ is a pair $(l, \pi_B^{\Omega})$, for some arrow
 $l:\Omega\times B \rightarrow \Omega $ in $\mathcal{E}$ satisfies the following conditions (as arrows in $\mathcal{E}$)\\
 $(1)$ $l\circ (l, \pi_B^{\Omega}) = l$;\\
$(2)$ $l\circ ( {\rm true}\times 1_B) = {\rm true} \circ  !_B$;\\
$(3)$ $l \circ\wedge_B =  \wedge\circ (l\circ (\pi_1, \pi_3), l\circ (\pi_2, \pi_3))$,
where $\pi_i$ is the $i$-th projection on $\Omega \times\Omega\times B,$ for $i=1, 2, 3$.
\end{lem}
By Lemma \ref{lop(e/b)}, for each topology $j$ on $\mathcal{E}$,
considering $l=j\circ \pi_{\Omega}^B$, it is easily seen that  $
j\times 1_B = (l, \pi_B^{\Omega})$ is a topology on
$\mathcal{E}/B$ which we  denote it by  $j_B$. In this case $j_B$
is called the {\it induced topology} on $\mathcal{E}/B$ by $j.$

One can simply see that if an arrow $k$ is a monomorphism in
$\mathcal{E}/B$, then $k$ as an arrow in $\mathcal{E},$ is too.
Also, for each monomorphism $k:f\rightarrowtail g$ in
$\mathcal{E}/B$, where $f:X\rightarrow B$ and $g : Y\rightarrow B$
in $\mathcal{E}$, we can observe
\begin{equation}\label{closure}
\widetilde{f} \stackrel{\widetilde{k}}{\rightarrowtail} g = (\overline{X}\stackrel{g\overline{k}}{\longrightarrow} B)\stackrel{\overline{k}}{\rightarrowtail} g,
\end{equation}
where $\overline{(\,\cdot\,)}$ and $\widetilde{(\,\cdot\,)}$ are
the universal closure operators corresponding to $j$ and $j_B$ on topoi $\mathcal{E}$ and $\mathcal{E}/B$,
respectively, in which  whole and the middle   squares of the
following diagram are pullbacks in $\mathcal{E}$,
\begin{equation}\label{pul clos}
\SelectTips{cm}{}\xymatrix{\overline{X}\ar@{>->}[dd]_{\overline{k}}\ar[rrr]&&&1\ar[dd]^{{\rm true}}\\
&X\ar[r] \ar@{>->}[d]_{k} & 1\ar[d]^{{\rm true}}&\\
Y\ar@{=}[r]&Y \ar[r]_{{\rm char}(k)} & \Omega\ar[r]^{j}&\Omega}
\end{equation}
(for more details, see \cite{maclane}).  One can construct
$\widetilde{k}$ in $\mathcal{E}/B$, similar to the above diagram.

  Here, we proceed to improve \cite[Vol. III, Proposition 9.2.5]{handbook 1} as follows:
\begin{lem}\label{pi}
Let $j$ be a topology in a topos $\mathcal{E}$. For every object
$B$ of $\mathcal{E}$, the pullback functor $\Pi_B :
\mathcal{E}\rightarrow \mathcal{E}/B$ preserves and reflects:
denseness (closeness) and $j$-separated objects ($j$-sheaves).
\end{lem}
\noindent{\bf  Proof.} Let $j$ be a topology on $\mathcal{E}$ and
$B$ an object of $\mathcal{E}$. Preserving dense (closed)
monomorphisms and sheaves (separated objects) in $\mathcal{E}$ by
the pullback functor $\Pi_B$, is standard and may be found  in
\cite[Vol. III, Proposition 9.2.5]{handbook 1}. To prove the rest
of lemma, here we just show that $\Pi_B$ reflects dense (closed)
monomorphisms. To verify this claim, let  $g:A\rightarrow C$ be an
arrow in $\mathcal{E}$ for which  $\Pi_B(g)$ is a $j_B$-dense
($j_B$-closed) monomorphism. We show that $g$ is $j$-dense
($j$-closed) monomorphism. As $\Pi_B (g) = g\times \id_B$ being
monomorphism in $\mathcal{E}/B$, the arrow $g$ is monomorphism in
$\mathcal{E}$ as well. For, let  $f, h$ in $\mathcal{E}$ be two
arrows such that $gf= gh$, we will have
$$\begin{array}{rcl}
gf= gh &\Longrightarrow & (g\times \id_B)(f\times \id_B) = (g\times \id_B)(h\times \id_B)\\
&\Longrightarrow & f\times \id_B = h\times \id_B \quad (g\times \id_B ~{\rm is ~ a~ monomorphism})\\
&\Longrightarrow & f=h.
\end{array}$$Considering $\overline{(\,\cdot\,)}$ and $\widetilde{(\,\cdot\,)}$ as the universal closure operators corresponding to $j$ and $j_B$,
 respectively. We get
$$\begin{array}{rclr}
\widetilde{\Pi_B (g)}&=& \widetilde{g\times \id_B}&\\
 &= & \overline{g\times \id_B} & \quad ~~~~~~~~~({\rm by}\quad (\ref{closure}))\\
&= & \overline{g}\times \id_B, &
\end{array}$$
where the last equality is true since we have $g\times \id_B =
({\pi_C^B})^{-1}(g)$, and because of stability of universal
closure operators
 under pullbacks we get $\overline{({\pi_C^B})^{-1}(g)} = ({\pi_C^B})^{-1}(\overline{g})$. The above equalities imply that if  $\Pi_B(g)$ is $j_B$-dense
 ($j_B$-closed) monomorphism in $\mathcal{E}/B$, then $g$ is $j$-dense ($j$-closed) monomorphism in $\mathcal{E}$.  $\qquad\square$

For any   topology $j$  on a topos $\mathcal{E}$, consider
$\mathcal{M}$ as the class of all $j$-dense  monomorphisms in
$\mathcal{E}$. Also, we denote  by $\mathcal{M}^{\perp}$ the class
of all arrows $g:C\rightarrow D$  in $\mathcal{E}$ such that for
any $f:A\rightarrow E$ in $\mathcal{M}$ and every commutative
square as in
\begin{equation}\label{factoriz}
\xymatrix{
A \ar[r] ^{u} \ar[d]_{f}& C \ar[d]^{g} \\E \ar[r]_{v} \ar@{-->}[ur]|w & D }
\end{equation}
there exists a unique arrow $w : E\rightarrow C$ in (\ref{factoriz}) such that the resulting triangles are commutative. In this case, we say that
 $g$ is {\it right orthogonal}  to $f.$ Moreover, we say that the pair $(\mathcal{M} , \mathcal{M}^{\perp})$ forms a {\it factorization system}
  in $\mathcal{E}$ if any arrow $f$ in  $\mathcal{E}$ factors as $f= me$, where $m\in \mathcal{M}$ and $e\in \mathcal{M}^{\perp}$ (for more information,
  see~\cite{Adamek}).
\begin{lem}\label{assum}
Let $j$ be a topology on a topos $\mathcal{E}$.  Then for each object $B$ of $\mathcal{E}$, we have $\mathcal{M}_B^{\perp}\subseteq \mathcal{M}^{\perp},$
where $\mathcal{M}_B$ is the class of all $j_B$-dense  monomorphisms in $\mathcal{E}/B$.
\end{lem}
{\bf  Proof.}  By Lemma \ref{pi} we get $\mathcal{M}_B\subseteq
\mathcal{M}$. To reach the conclusion, let $h:f\rightarrow g$ be
an arrow in $\mathcal{M}_B^{\perp}$,
  where $f:D\rightarrow B$ and $g:E\rightarrow B$ are arrows in $\mathcal{E}$. Now, consider the commutative square
\begin{eqnarray}\label{al}
\xymatrix{A \ar[r] ^{u} \ar[d]_{m}& D \ar[d]^{h} \\ C \ar[r]_{v} & E }
\end{eqnarray}
where $m:A\rightarrow C$ is in  $\mathcal{M}$. Since by  Lemma \ref{pi} the arrow $m:fu\rightarrow gv$   in $\mathcal{E}/B$ belongs to $ \mathcal{M}_B$
 and $h\in \mathcal{M}_B^{\perp}$,  there exists a unique arrow $w:gv\rightarrow f$ in $\mathcal{E}/B$ such that the following diagram commutes
\begin{eqnarray}\label{alii}
\xymatrix{fu \ar[r] ^{u} \ar[d]_{m}& f \ar[d]^{h} \\ gv \ar[r]_{v}\ar@{-->}[ur]|w & g }
\end{eqnarray}
The arrow $w:C\rightarrow D$ (as an arrow in $\mathcal{E}$) which
commutes the resulting triangulares, is unique in the diagram
(\ref{al}). To prove this, let $k : C\rightarrow D$ be an arrow in
$\mathcal{E}$ such that $km = u$ and $hk = v.$  Now,  we have $ fk
= (gh)k = gv$, so $k:gv\rightarrow f$ is an arrow in
$\mathcal{E}/B$ making all triangles in (\ref{alii}) commutative.
Thus, $k= w$ and  the proof is complete. $\qquad\square$
\begin{df}\label{enough c-dense}
{\rm Let $j$ be a topology on a topos $\mathcal{E}$. We say that  $\mathcal{E}$ {\it has  enough  $j$-sheaves} if for every object $A$ of  $\mathcal{E}$
there is a $j$-dense monomorphism $A\rightarrowtail F$ where $F$ is a $j$-sheaf.}
\end{df}

Following \cite{Adamek} a class $\mathcal{M}$ of morphisms in
$\mathcal{E}$ is  a {\it  left cancelable class}  if $g f\in
\mathcal{M}$ implies $f\in \mathcal{M}$. In the following, we
summarize the relation between  left cancelable  $j$-dense
monomorphisms  and factorization systems in a topos $\mathcal{E}$
and its slices.
\begin{thm}\label{assm}
Let  $j$ be a   topology on a topos $\mathcal{E}$. Assume that   for any object $B$  of $\mathcal{E}$, the class $\mathcal{M}_B$ of all $j_B$-dense
 monomorphisms in $\mathcal{E}/B$ be  left cancelable. Then  the following are equivalent:\\
{\rm (i)}  for any object $B$  of $\mathcal{E}$, $(\mathcal{M}_B , \mathcal{M}^{\perp}_B)$ is a factorization system in $\mathcal{E}/B$;\\
{\rm (ii)} for any object $B$  of $\mathcal{E}$, $\mathcal{E}/B$   has  enough  $j_B$-sheaves;\\
{\rm (iii)} for any object $B$  of $\mathcal{E}$, any object  of $\mathcal{E}/B$ is $j_B$-separated;\\
{\rm (iv)}  for any object $B$  of $\mathcal{E}$, any object  of $\mathcal{E}/B$ is $j_B$-sheaf;\\
{\rm (v)}   any object  of $\mathcal{E}$ is $j$-sheaf;\\
{\rm (vi)} any object  of $\mathcal{E}$ is $j$-separated;\\
{\rm (vii)} $\mathcal{E}$  has  enough  $j$-sheaves;\\
{\rm (viii)}  $(\mathcal{M} , \mathcal{M}^{\perp})$ is a
factorization system in $\mathcal{E}$.
\end{thm}
{\bf  Proof.}  That any $j$-sheaf is $j$-separated in
$\mathcal{E}$  yields that (v) $\Longrightarrow$ (vi)  holds.

 (vi) $\Longrightarrow$ (v). That  any object  of $\mathcal{E}$ is $j$-separated it follows that ${\bf Sep}_j(\mathcal{E})$
  is the topos $\mathcal{E}$ and then, every $j$-separated object is a $j$-sheaf as in  \cite[Theorem 2.1]{jhonstone1}.

(iii) $\Longrightarrow$ (vi). Setting $B=1$, then  any object  of $\mathcal{E}$ is $j$-separated.

(vi) $\Longrightarrow$ (iii). The claim follows immediately from
the fact that for any object $B$  of $\mathcal{E}$,
\begin{equation*}
{\bf Sep}_{j_B}(\mathcal{E}/B) \cong {\bf Sep}_j(\mathcal{E})/B.
\end{equation*}
 (see also \cite{jhonstone}).

(viii) $\Longrightarrow$ (vii).  By  (viii), for any object $A$ of
$\mathcal{E}$, the unique arrow $!_A:A\rightarrow 1$ factors as
$$\xymatrix{A~\ar[rr]^{!_A}\ar@{>->}[dr]_{m} &&1\\ &C\ar[ur]_{!_C}}$$
where $!_C\in  \mathcal{M}_1^{\perp} = \mathcal{M}^{\perp}$ and
$m\in  \mathcal{M}_1 =\mathcal{M}$. We remark that it is easy to
check that for any object $B$ of $\mathcal{E}$, $j_B$-sheaves in
$\mathcal{E}/B$ are exactly the class of all objects of
$\mathcal{E}/B$ which  belong to $\mathcal{M}_B^{\perp}.$  Since
$!_C$ is an object in $\mathcal{E}/1 = \mathcal{E}$ which is in
$\mathcal{M}^{\perp}_1 $, so $!_C$ is a $j_1$-sheaf, or
equivalently, $C$  is a $j$-sheaf.

(vii) $\Longrightarrow$ (viii). Consider an arrow $f:A\rightarrow B$ in $\mathcal{E}$. By using (vii), there exists a  $j$-dense monomorphism
$\iota:A\rightarrowtail F$, where  $F$ is a $j$-sheaf in $\mathcal{E}$. Now, we factor $f$ as the composite  arrow
 $A \stackrel{(\iota, f)}{\longrightarrow} F\times B\stackrel{\pi_B^F}{\longrightarrow} B$. Since $\pi_F^B(\iota, f) = \iota\in
 \mathcal{M}$ and $\mathcal{M}$ is a left cancelable class, so $(\iota, f)\in \mathcal{M}$. Also, $F$ being $j$-sheaf,  by Lemma \ref{pi}
 we have $\pi_B^F$ is a $j_B$-sheaf in $\mathcal{E}/B.$ By Lemma  \ref{assum} we have $\pi_B^F\in \mathcal{M}_B^{\perp}\subseteq \mathcal{M}^{\perp}$,
  as required.

(vi) $\Longrightarrow$ (vii). First of all we know  that any
$j$-separated object  of $\mathcal{E}$ can be embedded into a
$j$-sheaf (see, e.g.~\cite[Proposition V.3.4]{maclane}).  Let $A$ be an object of $\mathcal{E}$. Then, by
assumption $A$ is $j$-separated, and there exists an embedding
$A\stackrel{\iota}{\rightarrowtail} F$, where $F$  is a $j$-sheaf.
Now, take the closure of $A$ in $F$. Since $\overline{A}$ is closed
in $F$, by \cite[Lemma V.2.4]{maclane}, it is a $j$-sheaf. Since $A$
is $j$-dense in $\overline{A}$  we get the result.

(vii) $\Longrightarrow$ (vi). By  assumption for any object $A$ of
$\mathcal{E}$, there is a $j$-dense monomorphism $A\rightarrowtail
F$ in $\mathcal{E},$ where $F$ is a $j$-sheaf. Since any subobject
of a $j$-sheaf is $j$-separated so $A$ is  $j$-separated.

For any object $B$  of $\mathcal{E}$, setting $\mathcal{E}/B$
instead of $\mathcal{E}$ in (v), (vi), (vii) and (viii), we drive
(i) $\Longleftrightarrow$ (ii) $\Longleftrightarrow$ (iii)
$\Longleftrightarrow$ (iv). $\qquad\square$

In the following, we will introduce two main classes of dense
monomorphisms in a topos $\mathcal{E}$.
\begin{rk}\label{rkk}
{\rm By diagram~(\ref{pul clos}), one can
easily obtain: \\
 {\rm (i)} Let  $j = \id_{\Omega}$ be the trivial topology  on
$\mathcal{E}$. Then $j$-dense monomorphisms are  only the identity
maps. Therefore, any object of $\mathcal{E}$ is a $j$-sheaf. Also,
$j$-closed monomorphisms are exactly all
monomorphisms.\\
{\rm (ii)}  Let  $j$ be the topology ${\rm true} \circ !_{\Omega}$
on $\mathcal{E}$, that is, the characteristic map of
$\id_{\Omega}$. Then, $j$-dense  monomorphisms are exactly all
monomorphisms. Furthermore, $j$-closed  monomorphisms are just the
identity maps.}
\end{rk}

Recall~\cite{Adamek} that  $(Mono , Mono^{\square})$  is a weak
factorization system in any topos $\mathcal{E}$, where $Mono$ is the class
 of all monomorphisms in  $\mathcal{E}$. By Remark~\ref{rkk}(ii), the class $Mono$ is the class of all $j$-dense monomorphisms with
 respect to the topology $j = {\rm true} \circ !_{\Omega}$ on
$\mathcal{E}.$ Since the class $Mono$ is   left cancelable, so we
can obtain a special case of Theorem~\ref{assm} as follows.
(Notice that   by
 Lemma \ref{pi}  for the topology
  $j= {\rm true} \circ !_{\Omega} $ and any object $B$  of $\mathcal{E}$,
  the class $Mono_B$ will be  all monomorphisms in  $\mathcal{E}/B$.)
\begin{cor}\label{tnm}
For the topology $j= {\rm true} \circ !_{\Omega}$ on a topos
$\mathcal{E}$,  the following are equivalent:\\
{\rm (i)}   for any object $B$  of $\mathcal{E}$, $(Mono_B , Mono^{\perp}_B)$ is a factorization system in $\mathcal{E}/B$; \\
{\rm (ii)}  for any object $B$  of $\mathcal{E}$, $\mathcal{E}/B$   has  enough  $j_B$-sheaves;\\
{\rm (iii)}  for any object $B$  of $\mathcal{E}$, any object  of $\mathcal{E}/B$ is $j_B$-sheaf;\\
{\rm (iv)} for any object $B$  of $\mathcal{E}$, any object  of $\mathcal{E}/B$ is $j_B$-separated;\\
{\rm (v)}   any object  of $\mathcal{E}$ is $j$-sheaf;\\
{\rm (vi)}  any object  of $\mathcal{E}$ is $j$-separated;\\
{\rm (vii)} $\mathcal{E}$  has  enough  $j$-sheaves;\\
{\rm (viii)}  $(Mono , Mono^{\perp})$ is a factorization system in
$\mathcal{E}$.
\end{cor}
\section{Sheaves and sections of an arrow }
In this section, among other things, we investigate  a relationship between sheaves and sections of an arrow in a topos $\mathcal{E}$.
We start to remind \cite{Cagliari} that   for any object $B$ of
$\mathcal{E}$, the pullback functor $\Pi_B:\mathcal{E}\rightarrow
\mathcal{E}/B$ has a right adjoint $S:\mathcal{E}/B\rightarrow
\mathcal{E}$ as for any $f:X\rightarrow B$ we have the following
pullback
\begin{equation}\label{S(f) dia.}
\xymatrix{S(f)\ar[r]  \ar[d]_{}& 1 \ar[d]^{i_B} \\ X^B \ar[r]_{f^B} & B^ B }
\end{equation}
where $i_B$ is the transpose of $\id_B:1\times B\cong B\rightarrow B$ and $f^B$ is the transpose of the composition arrow
$X^B\times B \stackrel{ev_X}\longrightarrow X\stackrel{f}\longrightarrow B$ by the exponential adjunction $(-)\times B\dashv (-)^B$;
 that is, $ev_B (i_B\times \id_B ) = \id_B$ and $ev_B(f^B\times \id_B ) = f ev_X$, where the natural transformation
 $ev:(-)^B\times B\rightarrow (-)$ is the counit of  the exponential adjunction.
In fact, in the
Mitchell-B$\acute{{\rm e}}$nabou language, we
can write
$$S(f)=\{h \ | \ (\forall c\in B) \  f\circ (h(c)) = c \}.$$
This means that we can call $S(f)$  {\it the object of sections} of
$f$.

Since any retract of an object in a topos (or in an arbitrary category) is an equalizer, so the
topos ${\bf Sh}_j(\mathcal{E})$ is closed under retracts.
Furthermore, as $\Pi_B \dashv S,$ by Lemma~\ref{pi} we have that the
pullback functor $\Pi_B$ preserves dense monomorphisms, so $S$
preserves sheaves (for details, see~\cite[Corollary 4.3.12]{jhonstone}). (Roughly,  for any object $B\in \mathcal{E}$ and any
 adjoint $F\dashv G : \mathcal{E}\rightarrow \mathcal{E}/B$ one can easily checked that the functor $G$ preserves sheaves whenever $F$
 preserves  dense monomorphisms.)

In the following theorem we will find a relationship between sheaves in $\mathcal{E}/B$ and the object of sections of an arrow.
\begin{thm}\label{sectins}
Let $j$ be a topology  on a topos $\mathcal{E}$ and  $f:X\rightarrow
B$ be an object  of  $\mathcal{E}/B$. Then, $f$ is a $j_B$-sheaf in
$\mathcal{E}/B$, whenever the graph of $f$ which stands for the monomorphism $(\id_X,
f):f\rightarrowtail \pi_B^X$ in  $\mathcal{E}/B$, is a section  as well as
$S(f)$ is a $j$-sheaf in  $\mathcal{E}$.
\end{thm}
\noindent{\bf  Proof.} We recall that in \cite{Cagliari} it was
proved  if $(\id_X, f)$ is a section in  $\mathcal{E}/B$, then $f$ is
a retract of $\pi_B^{S(f)}$  in  $\mathcal{E}/B$. As $S(f)$ is a
$j$-sheaf, by Lemma \ref{pi}, $\pi_B^{S(f)}$  is a $j_B$-sheaf in
$\mathcal{E}/B$. But  ${\bf Sh}_{j_B}(\mathcal{E}/B)$ being closed
under retracts, therefore $f$ is a $j_B$-sheaf in
$\mathcal{E}/B$. $\qquad\square$

To the converse of Theorem \ref{sectins}, that the section functor $S$ preserves sheaves it yields that if  $f:X\rightarrow B$ be a $j_B$-sheaf
 in  $\mathcal{E}/B$,
 then $S(f)$ is a
$j$-sheaf in  $\mathcal{E}.$
Also, by Remark~\ref{rkk}(ii),  for $j = {\rm
true} \circ !_{\Omega}$, the monomorphism $(\id_X, f) : f\rightarrowtail \pi_B^X$  is
$j_B$-dense in  $\mathcal{E}/B$ and then for a $j_B$-sheaf $f:X\rightarrow B$, it will be a section in  $\mathcal{E}/B$.

In the rest of this section, for a small category  $\mathcal{C}$ we restrict our attention to obtain a version of Theorem \ref{sectins} for
 injective presheaves in  trivial slices of the presheaf topos $\widehat{\mathcal{C}}
={\bf Sets}^{\mathcal{C}^{op}}$ which is close to the version over
$j$-sheaves for the topology $j = {\rm true} \circ !_{\Omega}$  on
$\widehat{\mathcal{C}}.$ (See Proposition~\ref{sei} below.) Note
that the topology $j = {\rm true} \circ !_{\Omega}$ on
$\widehat{\mathcal{C}}$ is associated to the
 {\it chaotic or indiscrete Grothendieck topology} on $\mathcal{C}$.\\
 Recall \cite{maclane} that  in the presheaf topos $\widehat{\mathcal{C}}
={\bf Sets}^{\mathcal{C}^{op}}$, the exponential object $G^F$ is
defined in each stage $C$ of $\mathcal{C}$ as $G^F(C) = {\rm
Hom}_{\widehat{\mathcal{C}}}(Y(C) \times F, G)$, where $Y$ is the
Yoneda embedding, that is
$$Y:\mathcal{C}\rightarrow \widehat{\mathcal{C}}; \quad Y(C) = {\rm Hom}_{\mathcal{C}}(-, C).$$
Now, for  an arrow $\alpha:G\rightarrow F$ consider the arrows
$i_F:1\rightarrow F^F$ and $\alpha^F:G^F\rightarrow F^F$ in
$\widehat{\mathcal{C}}$ as the transposes of $\id_F:1\times F\cong
F\rightarrow F$ and $\alpha \circ ev_G:G^F\times F\rightarrow F$,
respectively,  by the exponential adjunction. We can observe
\begin{equation}\label{pls}
\forall C\in \mathcal{C},~~ (i_F)_C:1(C) = \{*\} \longrightarrow F^F(C);\quad (i_F)_C (*) = \pi_F^{Y(C)}.
\end{equation}
Also, for any two objects $C, D$ of $\mathcal{C}$, any $\gamma$ in
$G^F(C)$ and any $(k, y)$ in $Y(C)(D)\times F(D)$ we have
\begin{equation}\label{p}
(\alpha_C^F(\gamma))_D(k, y) = \alpha_D(\gamma_D(k, y)).
\end{equation}

Remind that a presheaf $G$ has a (unique)  global section which
means that in each stage $C$ of $\mathcal{C}$ there is a (unique)
element $\theta_C\in G(C)$ in such a way that for any arrow
$k:D\rightarrow C$ in $\mathcal{C}$ we have
\begin{equation}\label{C3}
G(k)(\theta_C)=\theta_D.
\end{equation}  Here, we find a special case that the exponential object and the
object of sections in $\widehat{\mathcal{C}}$ are exactly similar to  {\bf
Sets}. First, we express some lemma  required to achieve the goal.
 \begin{lem}\label{C}
Let $j$ be the topology ${\rm true} \circ !_{\Omega}$  on $\widehat{\mathcal{C}}$. Then,  the following assertions hold:\\
 {\rm (i)} For any $j$-sheaf $G$ in  $\widehat{\mathcal{C}}$, $G$ has a unique  global section. More generally, any injective presheaf $G$ of
 $\widehat{\mathcal{C}}$ has a  global section.\\
 {\rm (ii)} For any family $\{G_{\lambda}\}_{\lambda\in \Lambda}$  in  $\widehat{\mathcal{C}},$ the presheaf $G= \prod_{\lambda\in \Lambda} G_{\lambda}$
  is a $j$-sheaf (injective) in $\widehat{\mathcal{C}}$ iff  for all $\lambda\in \Lambda,$ $G_{\lambda}$ is a $j$-sheaf  (injective) in $\widehat{\mathcal{C}}.$
\end{lem}
{\bf  Proof.} (i) Let  $G$ be a $j$-sheaf in  $\widehat{\mathcal{C}}$ and  consider the coproduct object $G\sqcup 1$ in $\widehat{\mathcal{C}}.$
 By Remark~\ref{rkk}(ii),  there exists a unique natural transformation $\eta:G\sqcup 1\rightarrow G$ in $\widehat{\mathcal{C}}$ such that the following
  diagram commutes (if $G$ being  injective, the arrow $\eta$ is not necessarily unique)
$$\SelectTips{cm}{}\xymatrix{G \ar@{>->}[d]_{\iota}\ar[r]^{\id_G} & G  \\ \ar@{-->}[ur]_{\eta}G\sqcup 1}$$
where $\iota: G\rightarrow G\sqcup 1$ is the injection arrow.
Now, we will denote $\eta_C(*)$ by an element $\theta_C$ in $G(C)$ in each stage $C$ of $\mathcal{C}.$
Since $\eta:G\sqcup 1\rightarrow G$ is natural, so for any arrow $k:D\rightarrow C$ in $\mathcal{C}$ the following square commutes
$$\xymatrix{(G\sqcup 1)(D)\ar[r]^{~~\eta_D}  & G(D)\\
(G\sqcup 1)(C)\ar[r]^{~~~~\eta_C}\ar[u]^{(G\sqcup 1)(k)} & G(C)\ar[u]_{G(k)} }$$
Then, we have
$$\begin{array}{rcl}
G(k)(\theta_C)&=& G(k)(\eta_C(*)) \\
 &= & \eta_D ((G\sqcup 1)(k)(*))  \\
&= & \eta_D(1(k)(*)) =\theta_D.
\end{array}$$
This is the required result.

(ii) {\it Necessity.} Let $G$ be a $j$-sheaf (injective) in $\widehat{\mathcal{C}}.$ For any $\lambda, \mu\in \Lambda,$ we define
 $\alpha^{\lambda\mu}:G_{\lambda}\rightarrow G_{\mu}$ such that in each stage $C$ of $\mathcal{C}$ and for each $x\in G_{\lambda}(C),$
 we have $\alpha^{\lambda\mu}_C(x) = \theta_{C}^\mu,$ where $\theta_{C}^\mu$ is the $\mu$-th component of $\theta_{C}$ corresponding to $G$ in (i).
  Now, we will show that for any $\lambda, \mu\in \Lambda,$ $\alpha^{\lambda\mu}$ is a natural transformation in  $\widehat{\mathcal{C}},$ that is
 for any arrow $k:D\rightarrow C$ in $\mathcal{C}$ the following diagram is commutative
$$\xymatrix{G_{\lambda}(D)\ar[r]^{\alpha^{\lambda\mu}_D}  & G_{\mu}(D)\\
G_{\lambda}(C)\ar[r]^{\alpha^{\lambda\mu}_C}\ar[u]^{G_{\lambda}(k)}
& G_{\mu}(C)\ar[u]_{G_{\mu}(k)} }$$ For, consider an element $x\in
G_{\lambda}(C)$ we get
$$\begin{array}{rclr}
G_{\mu}(k)(\alpha^{\lambda\mu}_C(x))&=& G_{\mu}(k)(\theta_{C}^\mu)& \\
 &= & \theta_{D}^\mu &\quad ~~~~~ ({\rm by ~ (\ref{C3})})\\
 & = &  \alpha^{\lambda\mu}_D(G_{\lambda}(k)(x)). &
\end{array}$$
Now, for any $\lambda\in \Lambda,$ consider the family
$\{\gamma_{\mu}:G_{\lambda}\rightarrow G_{\mu}\}_{{\mu}\in
\Lambda}$ in $\widehat{\mathcal{C}}$
 such that for each $\lambda\not = \mu\in \Lambda$ we have $\gamma_{\mu} = \alpha^{\lambda\mu}$ and $\gamma_{\lambda} = \id_{G_{\lambda}}.$ Since $G$ is
 the product $\prod_{\lambda\in \Lambda} G_{\lambda},$ so there is a unique natural transformation $\gamma:G_{\lambda}\rightarrow G $ such that
 $p_{\mu}\gamma = \gamma_{\mu} $ and $p_{\lambda}\gamma = \id_{G_{\lambda}},$ for all $\lambda, \mu\in \Lambda$ and the projections $p_\lambda.$ Thus,
  for any $\lambda\in \Lambda, G_\lambda$ is a retract of the $j$-sheaf (injective) $G$ and then, $G_\lambda$ is a $j$-sheaf (injective).

{\it Sufficiency.} By the universal property of the product
presheaf $G,$ the unique arrow in the definition of a sheaf is
easily follows.$\qquad\square$

 We recall \cite{maclane} that  in each stage $C$ of $\mathcal{C}$ the object $ \Omega(C)$ of $\widehat{\mathcal{C}}$ is the set of all sieves on $C.$
   Also, the arrow ${\rm true}_C:1(C) = \{*\}\rightarrow \Omega (C)$  assigns to $*$, the {\it maximal  sieve} $t(C)$ of $ \Omega(C),$ that is all arrows with
   codomain $C$ of $\mathcal{C}.$
\begin{rk}
{\rm Note that the topology $j = {\rm true} \circ !_{\Omega}$ on $\widehat{\mathcal{C}}$ is the unique topology on $\widehat{\mathcal{C}}$ that satisfies
 Lemma~\ref{C}. To show this, for a $j$-sheaf $G$ of $\widehat{\mathcal{C}},$ consider the injection  $\iota: G\rightarrow G\sqcup 1$ in $\widehat{\mathcal{C}}.$
   In each stage $C$ of $\mathcal{C}$ we have ${\rm char} (\iota)_C(*) = \emptyset .$  Now, let $j$ be a topology on $\widehat{\mathcal{C}}.$ If $\iota$ is
 $j$-dense monomorphism, then in each stage $C$ of $\mathcal{C}$ we have $j_C(\emptyset ) = t(C).$  Now, for any sieve $S\in \Omega(C)$
 by Definition~\ref{loc. ope.}
  we get
$$\begin{array}{rcl}
t(C)   &= & j_C(\emptyset) = j_C(\emptyset\cap  S) \\
 &=&  j_C(\emptyset) \cap j_C(S)  = t(C) \cap j_C(S) = j_C(S).
\end{array}$$
Thus, $j_C$ is the constant function on $t(C),$ as required.}
\end{rk}

Let $F$ be the constant presheaf on a set $A.$   One can easily checked that the exponential adjunction $(-)\times F\dashv (-)^F$ is determined by, for any
 presheaf $G$ in $\widehat{\mathcal{C}},$ the exponential presheaf $G^F$ assigns to any object $C$ of $\mathcal{C}$,  the hom-set
 ${\rm Hom}_{{\bf Sets}} (A, G (C))$ and to any arrow $f:C\rightarrow D$ of $\mathcal{C}$,  the function
 $$G^F(f) : {\rm Hom}_{{\bf Sets}} (A, G (D)) \longrightarrow {\rm Hom}_{{\bf Sets}} (A, G (C)) $$ given by $G^F(f) (g) = G(f) \circ g$.
 As any function  $f:A\rightarrow G(C)$ can be considered as a sequence $(x_a)_{a\in A}\in \prod_A G(C),$ it yields that one has
  \begin{equation}\label{C6}
\forall C\in \mathcal{C},\quad G^F(C)\cong \prod_A G(C).
\end{equation}
By~(\ref{C6}), (\ref{pls}) and~(\ref{p}), it is convenient to see
that for each arrow $\alpha:G\rightarrow F$ in
$\widehat{\mathcal{C}}$ in
 which $F$ stands for the constant presheaf on a set $A,$  we get
\begin{equation}\label{C5}
 \forall C\in \mathcal{C},\quad S(\alpha) (C) \cong \prod_{a\in A}  \alpha_C^{-1} (a).
\end{equation}

Now, we will extract a special case of Theorem
\ref{sectins} in $\widehat{\mathcal{C}}$. First, let $\alpha:G\rightarrow F$ be an arrow in  $\widehat{\mathcal{C}}$ in which $F$ is
the constant presheaf on a set $A.$ For each element $a$ of $A,$ consider the subpresheaf $H_a$  of $G$ such that $H_a(C)=\alpha_C^{-1} (a)$,
 for any object $C$ of $\mathcal{C}.$ Since limits in $\widehat{\mathcal{C}}$ are constructed
pointwise, so  (\ref{C5}) shows that $S(\alpha) \cong \prod_{a\in A} H_a.$
\begin{pro}\label{C9}
Let  $j$ be the topology  ${\rm true} \circ !_{\Omega}$ on $\widehat{\mathcal{C}}$ and $\alpha:G\rightarrow F$ an arrow in  $\widehat{\mathcal{C}}$,
 where $F$ is the constant presheaf on a set $A$. Then, $\alpha$ is a
$j_F$-sheaf  in $\widehat{\mathcal{C}}/F$  iff the monomorphism $(\id_G, \alpha):\alpha\rightarrowtail \pi_F^G$ is a section in  $\widehat{\mathcal{C}}/F$
as well as
for any $a\in A,$ the subpresheaf $H_a$ of $G$ is a $j$-sheaf in $\widehat{\mathcal{C}}.$
\end{pro}
{\bf  Proof.} We deduce the result by Theorem \ref{sectins}, Lemma~\ref{C}(ii) and (\ref{C5}). $\qquad\square$

Since in topoi regular monomorphisms are exactly monomorphisms, so by \cite[Theorem 1.2]{Cagliari}, Lemma~\ref{C}(ii) and (\ref{C5}),
the following now gives which we are interested in.
\begin{pro}\label{sei}
Let   $\alpha:G\rightarrow F$ be an arrow in  $\widehat{\mathcal{C}}$, where $F$ is the constant presheaf on a
set $A$.  Then, $\alpha$ is  injective in $\widehat{\mathcal{C}}/F$ iff the monomorphism $(\id_G,
\alpha):\alpha\rightarrowtail \pi_F^G$ is a section in  $\widehat{\mathcal{C}}/F$ as well as for any $a\in A,$ the subpresheaf $H_a$ of $G$ is injective.
\end{pro}

In the case when $\mathcal{C}$ is a monoid, we obtain
\begin{ex}
{\rm Let $M$ be a monoid and $M$-{\bf Sets} the topos of all (right) representations of a fixed monoid $M.$ Since  $M$ is a small category with just
 one object, for two $M$-sets $X, B$ we have $X^B = {\rm Hom}_{M-{\bf Sets}}(M\times B, X)$, where $M\times B$ has the componentwise action. Hence,
  by (\ref{pls}) and (\ref{p}), for any equivariant map $f:X\rightarrow B$, in the diagram (\ref{S(f) dia.}) we observe
\begin{equation}\label{qls}
i_B (*) = \pi_B^{M}:M\times B\rightarrow B,
\end{equation}
and
\begin{equation}\label{q}
\quad\forall h\in X^B, ~\forall (m, b)\in M\times B,~~ (f^B(h))(m, b) = fh (m,b).
\end{equation}
Note that one writes  any
equivariant map $h:M\times B\rightarrow X$ in $X^B$ as a sequence
$((x_{m,b})_{b\in B})_{m\in M}$, consisting of  elements $x_{m,b} =
h(m, b)$ of $X$,  for any $(m, b)\in M\times B$. Also, $h$ being
equivariant map means that
$$\forall n,m\in M, \forall b\in B, \quad x_{mn, bn} = x_{m,b}n.$$
Hence,   we obtain that $X^B$ is equal to
\begin{equation}\label{as}
  \{  ((x_{m,b})_{b})_{m}\in \prod_{m\in M}\prod_{b\in B}X~\mid
\forall n,m\in M, \forall b\in B, x_{mn, bn} = x_{m,b}n \}.
\end{equation}
Now,  by (\ref{S(f) dia.}), (\ref{qls}) and (\ref{q}) we have
$$\begin{array}{rcl}
 S(f) &=& \{ ((x_{m,b})_{b})_{m}\in X^B~~|~ f^B(((x_{m,b})_{b})_{m}) = \pi_B^M  = ((b)_{b})_m\}\\
&=& \{ ((x_{m,b})_{b})_{m}\in X^B~~|~ ((f(x_{m,b}))_{b})_{m}   = ((b)_{b})_m\}\\
&=& \{ ((x_{m,b})_{b})_{m}\in X^B~~|~ \forall m\in M, \forall b\in B, x_{m,b}\in f^{-1} (b)\},
\end{array}$$
Hence, by (\ref{as}) we interpret a simple form of  underlying set of the $M$-set $S(f)$ in the topos $M$-{\bf Sets} as follows
$$\{ ((x_{m,b})_{b})_{m}\in \prod_{m\in M}\prod_{b\in B}f^{-1} (b)~ | ~\forall n,m\in M, \forall b\in B, x_{mn, bn} = x_{m,b}n \}.$$

 If $B$ has the trivial action $\cdot$, that is $\cdot = \pi_1 : B\times M\rightarrow B$ the first projection, then by (\ref{C6}) and (\ref{C5})
 we can obtain $X^B\cong \prod_BX$ and $S(f) \cong \prod_{b\in B}  f^{-1} (b)$.

 Furthermore, recall \cite{maclane} that for a group  $G$ and two $G$-sets $X, B$, we have
\begin{equation}\label{qq}
X^B=\{h:B\rightarrow X|~~h~{\rm is~ a~ function}\}\cong \prod_{B}X
\end{equation}
as two sets. According to the action on $X^B$, under the isomorphism (\ref{qq}),
the action on $\prod_{B}X$ is given by $ (x_b)_{b\in B}\cdot g = (x_{bg^{-1}}\cdot g)_{b\in B}$, for any $g\in G$ and $(x_b)_{b\in
B}\in \prod_{B}X$. Also, by (\ref{qq}) for any equivariant map
$f:X\rightarrow B$ in $G$-{\bf Sets}, in a similar way to
(\ref{C5}), we have $S(f) \cong \prod_{b\in B} f^{-1}(b)$.}
\end{ex}
\section{$j$-essential extensions in a topos}
This section is devoted to introduce  a class of monomorphisms in
an elementary topos, which we call these `$j$-essential
monomorphisms'. We present some equivalent forms of these and some
their properties. Meanwhile, we prove that any presheaf in a
presheaf topos has a maximal essential extension.

Remind that a  monomorphism  $\iota:A\rightarrowtail B$  is called
 {\it essential } whenever for each arrow $g:B\rightarrow C$
such that $g\iota$ is a monomorphism, then $g$ is a monomorphism
also. Now, we define a $j$-essential monomorphism  in a topos
 $\mathcal{E}$ as follows.
\begin{df}\label{j-essential mono. def}
{\rm For a topology $j$ on $\mathcal{E}$, a monomorphism
$\iota:A\rightarrowtail B$ is called {\it $j$-essential } whenever
it is $j$-dense
 as well as essential. In this case, we say that $B$ is a {\it $j$-essential extension} of $A$ and we
write $A\subseteq _{j} B$.}
\end{df}

We shall say an arrow $f: A\rightarrow B$ in   $\mathcal{E}$ is {\it $j$-dense} whenever the subobject $f(A)$, which is the image of $f$, is $j$-dense
 in $B$. In this way, any epimorphism in $\mathcal{E}$ becomes
 $j$-dense. ( For the definition of image of an arrow in a topos, see~\cite{maclane}.)

The following gives some equivalent definitions of  $j$-essential  monomorphisms  in a topos $\mathcal{E}$.
\begin{lem}\label{j-essential mono. equi. def.}
Let  $j$ be a topology on $\mathcal{E}$ and $\iota:A\rightarrowtail B$ a $j$-dense monomorphism. Then, the following are equivalent:\\
{\rm (i)} for any $g:B\rightarrow C$, $g$ is a monomorphism whenever $g\iota$ is a monomorphism;\\
{\rm (ii)} for any $g:B\rightarrow C$, $g$ is a $j$-dense monomorphism whenever $g\iota$ is a $j$-dense monomorphism;\\
{\rm (iii)} for any $g:B\rightarrow C$, $g$ is a monomorphism
whenever $g\iota$ is a $j$-dense monomorphism.
\end{lem}
{\bf  Proof.} (i) $\Longrightarrow$ (ii) and (iii) $\Longrightarrow$ (ii) are proved by~\cite[A.4.5.11(iii)]{jhonstone}.\\ (ii) $\Longrightarrow$
(iii) is clear.\\
(ii) $\Longrightarrow$ (i). Consider an arrow $g : B\rightarrow C$ for which $g\iota$ is monomorphism. We show that  $g$ is monomorphism also.
 Assume that $B\stackrel{k}{\twoheadrightarrow} g(B)\stackrel{m}{\rightarrowtail} C$
 is the image factorization of the arrow $g$. Since $g\iota = m (k\iota)$ and $g\iota$ is monomorphism, it follows that the arrow $k\iota$ is a
 monomorphism. Meanwhile, we get
 $$\begin{array}{rcl}
g(B) & =& k(B) ~~~~~~~~~~(\textrm{as} \ k\ \textrm{is epic})\\
&=& k(\overline{A})~~~~~~~~~~(\textrm{as} \ \iota \ \textrm{is dense})\\
&\subseteq & \overline{k(A)}\\
&\subseteq & g(B).
\end{array}$$
 Therefore, $g(B) = \overline{k(A)} = \overline{k\iota(A)}$. It follows that the compound monomorphism $k\iota : A\rightarrowtail g(B)$ is dense
 monomorphism and by (ii), $k$ is also.  That  $k$ is monomorphism and so isomorphism, yields that  $g$ is  monomorphism.$\qquad\square$

We point out that the proof of (ii) $\Longrightarrow$ (i) of Lemma \ref{j-essential mono. equi. def.} shows that  any composite
$k\iota$, for an epic $k$ and  a dense monomorphism $\iota$, is dense.

The follwing shows that $j$-essential monomorphisms in
$\mathcal{E}$  are closed under  composition.
\begin{pro}\label{compo. of j-ess. exte.}
Let  $j$ be a topology on $\mathcal{E}$. For two subobjects
$A\stackrel{\iota}{\rightarrowtail}A'\stackrel{~\iota'}{\rightarrowtail}B$
in $\mathcal{E}$,
 then  $A\subseteq _{j} B$ iff $A\subseteq _{j} A'$ and $A'\subseteq _{j} B.$
\end{pro}
{\bf  Proof.} By~\cite[13, A.4.5.11(iii)]{jhonstone}, one has $\iota' \iota$ is $j$-dense iff $\iota'$ and $\iota$ are $j$-dense.

{\it Necessity.} First, by Lemma~\ref{j-essential mono. equi.
def.}(i), we show that $A\subseteq _{j} A'$. To do so, consider an
arrow $f' : A'\rightarrow C$ for which $f'\iota$ is a
monomorphism.  Now, by~\cite[Corollary IV. 10. 3]{maclane}, the
object $C$  can be embedded into an injective object $D $ as in
$C\stackrel{\nu}{\rightarrowtail} D$ and hence there is an arrow
$\widetilde{f'} : B\rightarrow D$ such that $ \widetilde{f'}\iota'
=\nu f'.$ Since $A\subseteq _{j} B$ and $\widetilde{f'}\iota'
\iota = \nu f' \iota$ is a monomorphism, we deduce that
$\widetilde{f'}$ is a monomorphism. As $\widetilde{f'}\iota' = \nu
f'$ it follows that $f'$ is a monomorphism.

To prove $A'\subseteq _{j} B,$ choose an arrow $f : B\rightarrow C$ for which $f\iota'$ is a monomorphism. Then, $f\iota' \iota$
is also a monomorphism. Now $A\subseteq _{j} B$ implies that $f$  is a monomorphism, as required.

{\it Sufficiency.} Let  $f : B\rightarrow C$ be an arrow in $\mathcal{E}$ such that $f\iota' \iota$
is   a monomorphism. Since $A\subseteq _{j} A'$ and $(f \iota')\iota = f\iota' \iota$  is a monomorphism, it concludes that $f\iota'$
is a monomorphism. Using $A'\subseteq _{j} B,$  we achieve that $f$  is a monomorphism and hence $A\subseteq _{j} B$.$\qquad\square$

In the following, we achieve another property of $j$-essential monomorphisms in  $\mathcal{E}.$
\begin{lem}\label{embedding}
Let  $j$ be a topology on $\mathcal{E}$. If $A\subseteq _{j} B$ and $A$ is embedded in a $j$-sheaf  $F$, then $B$ also is  embedded in $F.$
\end{lem}
{\bf  Proof.} Let $\iota:A\rightarrowtail B$ be a  $j$-essential monomorphism and  $m:A\rightarrowtail F$ an arbitrary embedding.
 Since $F$ is a $j$-sheaf, there exists a unique morphism $f : B\rightarrow  F$ making the diagram below commutative;
\[\SelectTips{cm}{}\xymatrix{ A\ar@{>->}[r]^m\ar@{>->}[d]_{\iota} & F \\  B \ar@{-->}[ur]_f}\]
As $A\subseteq _{j} B$ being $j$-essential, $f$  is an embedding, as
required. $\qquad\square$

By Remark~\ref{rkk}(ii),  essential monomorphisms in a topos $\mathcal{E}$ are exactly $j$-essential monomorphisms in  $\mathcal{E}$
with respect to the topology $j = {\rm true} \circ !_{\Omega}$ on $\mathcal{E}.$

 Now, we would like to prove that any presheaf in $\widehat{\mathcal{C}}$ has a maximal essential extension.
\begin{thm}
Any presheaf in $\widehat{\mathcal{C}}$ has a maximal essential extension.
\end{thm}
{\bf  Proof.} Let $F$ be a presheaf in $\widehat{\mathcal{C}}$ and
$G$ an injective  presheaf into which $F$ can be embedded. By
Lemma~\ref{embedding}, we can assume that both $F$ and all its
essential extensions are subpresheaves of $G$.  Consider $\sum$ as
the set of all essential extensions of $F$ which is a poset under
subpresheaf inclusion $\subseteq$. Since the arrow $\id_F$ is an
essential extension of $F$, it follows that $\sum$ is non-empty.
If   $$\ldots \subseteq F_i\subseteq \ldots,$$ $i\in I,$ is a
chain in $\sum,$ then the subpresheaf $H$ of $G$ given by $H (C) =
\bigcup_{i\in I} F_i(C)$ for any object $C$ in $\mathcal{C}$ is an
upper bound of this chain. Now we show that $H$ lies in $\sum,$
i.e., $H$ is an essential extension of $F.$ To achieve this, let
$\alpha : H \rightarrow K$ be an arrow in $\widehat{\mathcal{C}}$
such that the restriction arrow $\alpha|_{F}$ is a monomorphism.
We prove that  $\alpha$ is a monomorphism. To verify this claim,
we show that for any $C\in \widehat{\mathcal{C}},$ the function
$\alpha_C : \bigcup_{i\in I} F_i(C) \rightarrow K(C)$ is one to
one.  Take $a , b\in \bigcup_{i\in I} F_i(C),$ $a\not = b.$ Then
there is a $j\in I$ such that $a , b\in F_j(C).$ Denote $\alpha|_{
F_j}$ by $\alpha_j.$ Since $F_j$ is an essential extension of $F$
and $\alpha_j|_{F} = \alpha|_{ F},$ it implies that $\alpha_j$ is
a monomorphism. Now
$$\alpha_C(a) = (\alpha_j)_C(a)\not = (\alpha_j)_C(b) = \alpha_C(b).$$
Therefore, $\alpha$ is a  monomorphism. Thus, $H\in \sum.$ Now it
follows from Zorn's Lemma that there is a maximal element $M$ in
$\sum.$ Then, $M$ is a maximal essential extension of $F.$
$\qquad\square$

It is straightforward to see that any  essential extension of $B$
can be embedded in any injective extension of $B.$

For a topology $j$  on a topos $\mathcal{E},$ by a {\it $j$-injective object}  we mean an injective object with respect to the class of
all $j$-dense monomorphisms in $\mathcal{E}.$

The following shows that the $j$-injective  presheaves
($j$-sheaves) in $\widehat{\mathcal{C}}$  have no proper
$j$-essential extension.
\begin{pro}\label{non- pro. j-ess. exte.}
Let  $j$ be a topology on $\widehat{\mathcal{C}}$ and  $F$ a
$j$-injective  presheaf ($j$-sheaf) in $\widehat{\mathcal{C}}.$
Then, $F$ has
 no proper  $j$-essential extension.
\end{pro}
{\bf  Proof.} Suppose that $G$ is a  proper $j$-essential
extension of $F$ and so $F$ is a $j$-dense subpresheaf of $G$ and
$F\not = G.$ Thus there is an object  $C$ of $\mathcal{C}$ such
that $G(C)\not \subset F(C)$ and then, an $a\in G(C)$ such that
$a\not \in F(C).$ Since $F$ is $j$-injective ($j$-sheaf) implies
that there is an arrow $\alpha : G\rightarrow F$ for which
$\alpha|_{F} = \id_{F}.$ That $a\not \in F(C)$ and $\alpha_C(a)\in
F(C)$ follows that $a\not = \alpha_C(a).$ But
$\alpha_C(\alpha_C(a)) = \alpha_C(a).$ Then, $\alpha_C$ and so
$\alpha$ is not a monomorphism
 although  $\alpha|_{F} = \id_{F}$ is. This shows that $G$ is not a  proper $j$-essential extension of $F$ and it is a contradiction.$\qquad\square$

The following shows that  the pullback functor $\Pi_B$ reflects $j$-essential extensions.
\begin{pro}
Let $j$ be a topology in a topos $\mathcal{E}$. For every object $B\in \mathcal{E}$, the pullback functor $\Pi_B : \mathcal{E}\rightarrow \mathcal{E}/B$
 reflects $j$-essential monomorphisms.
\end{pro}
\noindent{\bf  Proof.}  Let   $f:A\rightarrow C$ be an  arrow in  $\mathcal{E}$ such that $\Pi_B(f)$ is a $j_B$-essential monomorphism. We show that $f$ is a
$j$-essential monomorphism. By Lemma \ref{pi}, $f$ is a $j$-dense monomorphism in $\mathcal{E}$. Let $g:C\rightarrow D$ be an arrow in $\mathcal{E}$
 such that $gf$ is a monomorphism. We show that $g$ is too.
Since  $gf$ is a monomorphism, the arrow $(g\times \id_B) \Pi_B(f) = (gf)\times \id_B$ is also a  monomorphism. As $\Pi_B(f)$ is  $j_B$-essential,
 so $g\times \id_B$  is a monomorphism. Then $g$ is a monomorphism. This is the required result.$\qquad\square$

Recall~\cite{Madanshekaf} that a {\it weak topology} on a topos $\mathcal{E}$ is a morphism $j : \Omega\rightarrow \Omega$ such that:\\
(i) $j\circ {\rm true} = {\rm true}$;\\
(ii) $j\circ \wedge \leq  \wedge\circ (j\times j)$,
in which $\leq$ stands for the internal order on $\Omega.$
 Meanwhile, a weak  topology $j$ on $\mathcal{E}$ is said to be  {\it productive} if
$j\circ \wedge = \wedge\circ (j\times j)$.

In what  follows, we review the whole paper for a weak topology $j$ on a topos $\mathcal{E}$ instead of a topology.
\begin{rk}
{\rm Similar to~\cite[A.4.5.11(ii)]{jhonstone}, one can easily
check that for a weak topology $j$ on $\mathcal{E}$  pushouts also
preserve dense monomorphisms. Hence, we can obtain a version of
Lemma \ref{unique retract} for a weak  topology $j$ on
$\mathcal{E}$ as well. One can observe that completely analogous
assertions to Lemmas~\ref{j-essential mono. equi.
def.},~\ref{embedding}  and \ref{pi}, Proposition~\ref{non- pro.
j-ess. exte.}
 and Theorem  \ref{sectins}, hold for  a weak  topology $j$ on $\mathcal{E}$. But, by~\cite{Madanshekaf}, in the proof of Theorem  \ref{assm},
 the part (vi) $\Longrightarrow$ (vii) is true for  a productive weak  topology $j$ on $\mathcal{E}$. The rest parts of this proof satisfies
 for weak  topologies.\\
Recall~\cite{Madanshekaf} that, for a weak topology $j$ on
$\mathcal{E},$ it is convenient to see that if the composite
subobject $mn$ is dense then so are $m$ and $n$. In contrast with
 topologies~\cite[A.4.5.11(iii)]{jhonstone}, the
converse is not necessarily true. Hence, the sufficiency part of
Proposition~\ref{compo. of j-ess. exte.} does  not necessarily
hold for a weak topology  $j$ on a topos
 $\mathcal{E}.$ The necessity part of this proposition satisfies for a weak topology  $j$ as well.}
\end{rk}

\end{document}